\documentclass[11pt]{amsart}
\usepackage{graphicx}
\usepackage{amsfonts}
\usepackage{amscd}
\usepackage{amssymb}
\usepackage{alltt}

\newcommand{\ring}[1]{\mathbb{#1}}

\newcommand{\calL}{{\mathcal L}}

\newcommand{\calK}{{\mathcal K}}
\newcommand{\lie}[1]{{\bf\mathfrak #1}}
\def\infnty{\infty}
\def\mot{{K_0^v(\op{Mot}_{\ring{Q},\bar{\ring{Q}}})_{\op{loc},\ring{Q}}}}

\def\op#1{{\operatorname{#1}}}

\title{Orbital Integrals are Motivic}
\author{Thomas C. Hales}

\date{December 5, 2002}

\begin{document}

\theoremstyle{plain}
\newtheorem{thm}{Theorem}
\newtheorem{lemma}[thm]{Lemma}
\newtheorem{cor}[thm]{Corollary}
\newtheorem{prop}[thm]{Proposition}
\newtheorem{cond}[thm]{Condition}

\theoremstyle{definition}
\newtheorem{rem}[thm]{Remark}
\newtheorem{defn}[thm]{Definition}
\newtheorem{ex}[thm]{Example}
\newtheorem{estimate}[thm]{Estimate}

\begin{abstract}  This article shows that under general conditions,
$p$-adic orbital integrals of definable functions are represented by
virtual Chow motives. This gives an explicit example of the philosophy
of Denef and Loeser, which predicts that all ``naturally occurring''
$p$-adic integrals are motivic.
\end{abstract}
\maketitle

\section{Introduction}%
\footnote{I grant this paper to the public domain. No rights are
reserved by the author.} Denef and Loeser have introduced a theory of
arithmetic motivic integration \cite{DL}.   In this theory, general
families of $p$-adic integrals can be calculated as the trace of a
Frobenius operator on virtual Chow motives.  This article shows that
orbital integrals fit nicely into the general framework of Denef and
Loeser.  It describes a large class of orbital integrals that can be
computed by a Frobenius operator on virtual Chow motives.  Moreover,
there is an effective procedure to compute the virtual Chow motive from
the data defining the orbital integral.  In this sense, this article
gives an algorithm to compute a large class of orbital integrals.

The idea of using double cosets to compute motivic orbital integrals is
taken from J. Gordon's recent thesis (and a suggestion of Julee Kim).
The thesis proves that under general conditions, the character values of
depth zero representations of a $p$-adic group can be represented as
virtual Chow motives \cite{G}. This note can be viewed as an extension
of the methods of that thesis.

\section{Double coset bounds}

Recall that Pas has defined a first order language that is based on the
theory of valued fields \cite{P}.  It is a three-sorted language in the
sense of \cite{E}. The models of the three sorts are a valued field, a
residue field, and the additive group of integers (the value group),
augmented by $+\infnty$.  The language has function symbols $\op{ord}$
and $\op{ac}$ that are interpreted as the valuation and angular
component maps on a valued field.  The valuation is a map from the
valued field to $\ring{Z}\cup\{+\infnty\}$. The angular component map is
a function from the valued field to its residue field.

Let $\calK$ be a set of models of Pas's language.  We do not make any
assumptions about the residual characteristic of the fields in $\calK$.
We do, however, assume that each field is complete and Henselian. In
fact, the only cases of interest to us are locally compact
nonarchimedean fields (in brief, $p$-adic fields).  For example, we
could take $\calK$ to be models corresponding to the set of $p$-adic
fields $\ring{Q}_p$, or the set of fields $\ring{F}_p[[t]]$, or the set
of all $p$-adic fields. We write $\theta^k$ for an interpretation of
$\theta$ in the model $k$.  If $S$ is a finite set of prime numbers, let
$\calK_S$ be the set of $k\in\calK$ for which the residual
characteristic of $k$ is {\it not\/} in $S$.

Let $\theta(m)$ be a formula in Pas's language. Assume that $\theta$ has
no free variables of the valued field sort or of the residue field sort.
Assume that its free variables of the additive sort are contained in
$m=(m_1,\ldots,m_\ell)$.

Suppose that there exists a finite set $S$ of prime numbers satisfying
the following condition:

\begin{cond}  \label{cond:bound}
For every $k\in\calK_S$,
    $$
    \{m\in\ring{Z}^\ell \ : \ \theta^k(m)\}
    $$
is a bounded subset of $\ring{Z}^\ell$.
\end{cond}

\begin{thm} \label{thm:bound}
Under the stated conditions, there exists a finite set $S'$ of prime
numbers and a bounded subset $C \subset \ring{Z}^\ell$ such that for
every $k\in\calK_{S'}$ we have
    $$
    \{(m_1,\ldots,m_\ell)\in\ring{Z}^\ell \ : \
        \theta^k(m_1,\ldots,m_\ell)\}
        \subset C.
    $$
\end{thm}

In other words, by throwing away finitely many primes, the bound on the
subset in Condition~\ref{cond:bound} can be made independent of the
model.

\begin{proof}
Apply Pas's quantifier elimination on quantifiers of valued field sort
in the formula $\theta$.    We obtain an equivalent formula that
contains no variables of valued field sort. This formula is uniform in
the sense that it is independent of the model $k\in\calK$.  Pas assumes
that the residual characteristic is zero.  We achieve equivalent results
by throwing out finitely many residual characteristics.

The remaining terms of valued field sort are constants that are
definable in Pas's language without quantifiers. Such a constant must be
a rational number. Equalities in the valued field $a = b$ can be
replaced with an equivalent statement with the angular component map:
$\op{ac}(a-b)=0$. The infinite valuation condition
$\op{ord}(a)=+\infnty$ can be replace with $\op{ac}(a)=0$. Excluding
finitely many primes, each rational number appearing in the formula can
be assumed to be $0$ or a unit; that is, it is an integer in each model.
If $a\ne 0$, the condition $\op{ord}(a)=m$ is then equivalent to the
formula $m=0$, and if $a=0$ it is false (because we have already treated
$\op{ord}(a)=+\infnty$).   Similarly, the formula $\op{ac}(a)=\xi$ is
equivalent to $a'=\xi$, where $a'$ is the element in the residue field
that is the reduction of the integer $a$. In this way, we eliminate all
terms of valued-field sort, all function symbols $\op{ac}$ and
$\op{ord}$.

Writing the formula in disjunctive normal form we find that $\theta$ is
equivalent to
    \begin{equation}\bigvee_i (\psi_i \wedge L_i),\label{eqn:vee}\end{equation}
where $\psi_i$ is a formula with no free variables and containing only
constants and terms of the residue field sort, and $L_i$ is a formula
with free variables limited to $m$ and whose constants and terms are of
the value group sort.  We may assume that $+\infnty$ has been eliminated
from the language, so that each $L_i$ is a conventional Presburger
formula \cite{E}.

We partition the finite set of indices $i$ into two sets $B$ and $\bar
B$, where $B$ is the set of indices for which $L_i$ is a bounded subset
of $\ring{Z}^\ell$ and $\bar B$ the set of indices for which $L_i$ is
unbounded.  We can find a set $C\subset \ring{Z}^\ell$ that contains
    $$\cup_{i\in B} \{m\ :\ L_i(m)\}.
    $$
Note that $C$ is independent of the model.

Let $S'$ be the union of $S$ with the set of primes that were excluded
by Pas's quantifier elimination, together with the primes that were
excluded to make the rational constants in the valued field all units
(or 0).

If $k\in\calK_{S'}$, then $\theta$ is equivalent to
Formula~\ref{eqn:vee} in the interpretation $k$.  By hypothesis, the set
of solutions is bounded, which implies that each $\psi_i^k$, for
$i\in\bar B$, is false. Hence all the solutions lie in $C$, as desired.
\end{proof}

\section{Some virtual sets}

We write
$$\theta(x_1,\ldots,x_\ell,m_1,\ldots,m_{\ell'},\xi_1,\ldots,\xi_{\ell''})$$
for a formula in Pas's language, where all free variables are among
those listed, and the variables are of the valued field sort $x_i$,
value group sort $m_i$, and residue field sort $\xi_i$.

If $\theta(x_1,\ldots,x_\ell)$ is a formula in Pas's language with no
free variables of residue field sort and no free variables of the value
field sort, then we have the formula
    $$\begin{array}{lll}
        &\exists m \ \forall y_1,\ldots,y_\ell,x_1,\ldots,x_\ell\\
        &\quad (\theta(x_1,\ldots,x_\ell) \wedge
        (\op{ord}(y_i-x_i)\ge m, \text{ for } i=1,\ldots,\ell)\\
        &\quad \Rightarrow \theta(y_1,\ldots,y_\ell).
        \end{array}$$
This formula does not contain any free variables.  By quantifier
elimination, this is equivalent to something of the shape of
Formula~\ref{eqn:vee}.  There are no free variables of the value group
sort.  By Presburger quantifier elimination, the formulas $L_i$ can be
replaced with ``true'' $1=1$ or ``false'' $0=1$.  Thus, we may assume
that the formula is a formula in the first order theory of the residue
field sort.  We say that $\theta(x_1,\ldots,x_\ell)$ is {\it stable} if
the formula in the residue field sort is true for all pseudo-finite
fields.\footnote{A field is pseudo-finite if it has a single extension
of each degree and if each absolutely irreducible variety over the field
has a rational point.} If the formula $\theta$ is stable, then by
avoiding finitely many residue field characteristics, interpretations of
$\{x:\theta(x)\}$ are stable definable sets. ({\it Stable\/} is meant in
the sense of Denef and Loeser's motivic integration \cite{DL}.)

A {\it virtual set\/} is a class construct
    $$\{x : \theta(x)\},$$
where $\theta$ is a formula in Pas's language.  See \cite{VTF}, where
the corresponding notion for formulas in the first-order theory of rings
is discussed.

Let $\lie{g}$ be a virtual reductive Lie algebra in the sense of
\cite{VTF}; that is, a virtual set in Pas's language whose models are
reductive Lie algebras in the traditional sense.  Assume that $\lie{g}$
is the Lie algebra of a virtual reductive Lie group $G$.  (Again, this
is to be interpreted as a virtual set whose models are reductive
algebraic groups.)  For example, each split reductive group over
$\ring{Q}$ gives a virtual reductive group in Pas's language.

We define the space of {\it locally constant virtual functions\/} on
$\lie{g}$ to be the $\ring{Q}$-linear combinations of characteristic
functions of stable virtual sets in $\lie{g}$.   More precisely, take
the Grothendieck group generated by stable formulas in Pas's language,
and tensor with $\ring{Q}$.   Excluding finitely many primes, there is a
map from this ring to the ring of locally constant $\ring{Q}$-valued
functions on $\lie{g}^k$ for $p$-adic field models $k$.

Let $\ell$ be the rank of $\lie{g}$.  Let
$\op{span}(Y,e_1,\ldots,e_\ell)$ be the formula in the first-order
language of rings asserting that $Y$ lies in the span of the $e_i$.  The
set of regular elements $\lie{g}^{reg}$ is defined as a virtual set of
elements in $\lie{g}$
    $$
    \{ X \ : \ \exists e_1,\ldots,e_\ell\ \forall Y.\quad
        [X,Y]=0\ \Rightarrow\ \op{span}(Y,e_1,\ldots,e_\ell)\}.
    $$
(The quantifiers in this formula have been extended to range over
elements of $\lie{g}$, as was done in \cite{VTF}.)

The set of nilpotent elements in $\lie{g}$ is a subvariety, hence a
virtual subset. A regular semisimple element is defined as a regular
element that has $0$ as the only nilpotent element in the centralizer.
This condition is given by a formula, hence regular semisimple elements
form a virtual set $\lie{g}^{reg,ss}$.

In a split reductive algebra over $\ring{Q}$, we have a finite list of
definable subsets of $\lie{g}$ giving all of the proper parabolic
subgroups $P$ up to conjugacy.  The virtual set of regular semisimple
elliptic elements is given by
    $$
    \cap_P
    \{ X \in \lie{g}^{reg,ss}\ :\
    \text{There does not exist } g \in G \text{ such that } \op{Ad}\,g (X) \in
    P\}.
    $$

\begin{ex}  Consider the virtual set  in $GL(2)$:
    $$
    Z = \{\begin{pmatrix}a & b\\ u b &
    a\end{pmatrix}:\ \op{ord}(u)=0 \wedge \op{ord}(a)\ge0\wedge
    \op{ord}(b)\ge0\wedge a^2- u b^2\ne0\wedge b\ne0\}.
    $$
This is a virtual subset of the set of regular semisimple elements.  The
condition for ellipticity for $X\in Z$ is
    $$
    \exists A.  \ A Z - Z A \in \begin{pmatrix} *&*\\ 0&*\end{pmatrix}.
    $$
This is equivalent (excluding the prime $2$) to the virtual subset of
$Z$ given by the condition
    $$
    \exists \xi. \ \op{ac}(u)=\xi^2.
    $$
\end{ex}

\section{An integration formula}
\label{sec:integrate}

Let $G = \sqcup_t K t K$ be the Cartan decomposition of the set $G$ of
$p$-adic points of a split reductive group over a $p$-adic field. Assume
$t\in T$ the set of points of a split torus.   The method of Cartier
gives the following integration formula for functions supported on a
single double coset $K t K$.

\begin{thm}    Let $dk$ be a Haar measure on $K$.
Let $dg$ be a Haar measure on $G$.  Then there exists a constant
    $c(K t K)$ such that for all $f\in C_c(K t K)$.
    $$
    c(K t K) \int_{K\times K} f(k_1 t k_2) dk_1\,dk_2
        = \int_G f(g) dg.
    $$
\end{thm}

\begin{proof}  Follow Cartier \cite{C}.
\end{proof}

The constant can be obtained by taking $f$ to be the characteristic
function of $K t K$.  We find that
    $$
    c(K t K) \op{vol}(K,dk)^2= \op{vol}(K t K,dg).$$
Pick $dg$  so that its restriction to $K$ is $dk$, and normalize $dk$ so
that
    $$\op{vol}(K,dk) = |G(\ring{F}_q)| q^{-\op{dim}(G)}.$$
Let $[G]_q$ be the constant defined by the right-hand side of this
equation.

\section{The Main Result}

Let $\lie{g}$ be a split reductive group over $\ring{Q}$.  Let $G$ be a
split reductive group over $\ring{Q}$ with Lie algebra $\lie{g}$.  We
take $\lie{g}$ and $G$ to be represented explicitly as matrices in
$\lie{gl}(n)$ and $GL(n)$, respectively.  Let $O_v$ be the ring of
integers of a completion $L_v$ of a number field $L/\ring{Q}$.  We may
assume that $G(O_v) = GL(n,O_v)\cap G(L_v)$ is a hyperspecial maximal
compact subgroup for almost all places $v$ of $L$.

View $\lie{g}$ and $G$ as virtual subsets of the virtual sets
$\lie{gl}(n)$ and $GL(n)$.  Let $\lie{gl}(n,O)$ be the virtual subset of
$\lie{gl}(n)$ given by
    $$\{ (x_{ij})\ : \ \op{ord}(x_{ij})\ge0, \ \forall i,j\}.$$
Set $\lie{g}(O) = \lie{gl}(n,O)\cap \lie{g}$.

Let $f$ be a locally constant virtual function on $\lie{g}$.  Let $E$ be
a stable virtual set of regular semisimple elliptic elements. (Recall
that `stable' is used here in the sense of stability in motivic
integration, and not in the sense of stable conjugacy in group theory.)
Assume that $E\subset \lie{g}(O)$.

If $L$ is a finite extension of $\ring{Q}$, and $v$ is a place of $L$,
the completion $L_v$ at $v$ is (the domain of) a model of Pas's
language. Let $O_v$ be the ring of integers of $L_v$.  We write $f_v$
for the locally constant function in that model.  We write $E_v$ for the
definable subset of $\lie{g}(O_v)$ obtained from the virtual set $E$.
Assume that for almost all $L_v$, $E_v$ is a compact set. (By `almost
all' here and in what follows,  we mean all but finitely many residue
field characteristics for $L_v$.)

A ring $\mot$ of virtual Chow motives is constructed by Denef and Loeser
in \cite{DL}.  The ring $\mot$ is defined as the image of a ring
$K_0(\op{Sch}_{\ring{Q}})_{\op{loc}}$ in the Grothendieck ring of the
category of Chow motives over $\ring{Q}$.  The ring
$K_0(\op{Sch}_{\ring{Q}})_{\op{loc}}$ is the Grothendieck ring of Chow
motives over $\ring{Q}$, localized by inverting the class $\ring{L} =
[\ring{A}^1_{\ring{Q}}]$ of the affine line.  See \cite{DL} for details.
In this paper, a {\it virtual Chow motive} $M$ will always mean an
element of this ring.

If $v$ is a place of $L$, then there is a Frobenius operator
$\op{Frob}_v$ that acts on the $\ell$-adic cohomology groups of $M$. See
\cite[3.3]{DL}.   We write the alternating trace of Frobenius on the
cohomology as
    $$\op{trace}(\op{Frob}_v,M)\in \bar {\ring{Q}}_\ell.$$

We fix measures on $G_v$ (the set of $p$-adic points of $G$), and
$\lie{g}_v$ (the set of $p$-adic points of $\lie{g}$) by fixing
differential forms of top degree defined over $\ring{Q}$ on $G$ and
$\lie{g}$.  Normalize $dX$ so that $\op{vol}(\lie{g}(O_v))=1$. Normalize
$dg$ to restrict to $dk$ on $G(O_v)$.  Normalize $dk$ as in
Section~\ref{sec:integrate}.

\begin{thm}  For each $f$ and $E$,
there is a virtual Chow motive $M(f,E)$ over $\ring{Q}$ and a finite set
$S$ of primes with the following property with the following property:
If $L$ is a finite extension of $\ring{Q}$ and $v$ is a place of $L$
that does not lie over any $p\in S$, then
    $$
    \op{vol}(K,dk)
    \int_{E_v} \int_{G_v} f_v (g^{-1} X_v g) dg \,d X
        = \op{trace}(\op{Frob}_v, M(f,E)).
        $$
\end{thm}

\begin{rem} In general, a single elliptic element $X$ in a $p$-adic
field is not definable and the conjugacy class $O(X)$ is not definable.
Thus, it is not reasonable to ask for a virtual Chow motive to represent
the integral, except when averaged over $E$.
\end{rem}

\begin{rem} The same conclusion holds if we take our models to be
$\ring{F}_q((t))$ instead of $L_v$.  The Chow motive $M(f,E)$ is the
same in both cases.  Hence, we find as a corollary that the orbital
integrals are ``the same'' in zero characteristic and positive
characteristic.
\end{rem}

\begin{proof}  In the course of the proof we increase the size of the
finite set $S$ of primes several times without a change in notation.

For any given $L_v$, the integral over $G_v$ is a sum of integrals over
$K t K$, as we run over all double cosets.  For any given $X\in E_v$,
there are only finitely many double cosets involved (because $X$ is
elliptic). The collection of double cosets required is locally constant
in $X$, so there is a neighborhood of $X$ that on which the finite set
of cosets is the same as for $X$.  By compactness on $E_v$, there is a
finite cover by such neighborhoods.  It follows that there is a finite
set of double cosets depending only on $E_v$ (and $f_v$).

For simplicity, we may assume that $f_v$ is a characteristic function at
$L_v$ of a virtual set $D$.  Let $m=(m_1,\ldots,m_\ell)$ be a tuple of
integers with $\ell$ equal to the reductive rank. Let $T$ be a split
torus that gives at each $p$-adic place the split torus that appears in
the Cartan decomposition.  Let $\phi:\ring{G}_m^{\times\,\ell}\to T$ be
an isomorphism (defined over $\ring{Q}$) between a product of
multiplicative groups and the diagonal torus $T$.  Extend $\op{ord}$ to
$\op{ord}:T\to \ring{Z}^\ell$ by
    $$T\to \ring{G}_m^\ell \to^{\op{ord}^\ell} \ring{Z}^\ell.$$

Let $\psi(m,k_1,k_2,X)$ be the formula in Pas's language given by
    $$\begin{array}{lll}
    &\exists t.\quad X\in E\ \wedge\
        \op{Ad}(k_1 t k_2) X \in D\\
        &\quad\wedge\
            \op{ord}_v (t) = m\\
        &\quad\wedge\
            k_1,k_2\in \lie{g}(O).
    \end{array}
    $$
Being in Pas's language, the formula is independent of the model
$(L,v)$. Let $\psi(m)$ be the formula
    $$\exists\, k_1\, k_2\, X. \ \psi(m,k_1,k_2,X).$$

The formula $\psi(m)$ has a bounded number of solutions in each model
$\theta_v$. By Theorem~\ref{thm:bound}, there is a finite set
$C\subset\ring{Z}^\ell$ that contains the solutions in $m$ for almost
all models $L_v$.  This implies that we can pick a finite collection of
double cosets that works simultaneously for almost all models $L_v$. By
the Cartan decomposition, we may assume that $m_1\ge m_2\ge \cdots$ for
$m = (m_1,\ldots,m_\ell)\in C$.  Write $K_m$ for the double coset
    $$K\phi(\pi^m)K.$$

By the Cartier-style integration formula \ref{sec:integrate}, the
integral of the theorem can be written as
    $$
    \sum_{m\in C} [K_m:K] \int_{E_v\times K\times K}
        \op{char}(\psi(m,k_1,k_2,X))
    dk_1\,dk_2\,dX.
    $$

By the main theorem of Denef and Loeser \cite{DL}, for each $m$ in this
finite sum, there is a motive $M_m$ such that the trace of Frobenius at
$v$ computes this volume.  In their theorem, the $p$-adic volume is to
be computed with respect to the Serre-Oesterl\'e measure on the group
    $$G(O_v)\times G(O_v)\times \lie{g}(O_v).$$
This measure is invariant under analytic isomorphisms of the set. In
particular, it is invariant under the group action of
    $$G(O_v)\times G(O_v)\times \lie{g}(O_v)$$
on itself.  Hence it is a Haar measure.

To work out the normalization of the Haar measure, we note that the
Serrre-Oesterl\'e measure coincides with the counting measure modulo a
uniformizer, scaled by $q^{-\op{dim}(X)}$.  Thus, the volume of this
group is
    $$[G]_q\times[G]_q.$$
For the normalization of $dk$ and $dX$ chosen in
Section~\ref{sec:integrate}, the measure is also $[G]_q\times[G]_q$.
Hence the Serre-Oesterl\'e measure coincides with the chosen Haar
measure on this group.

The constant $[K_m:K]$, the number of cosets, is a polynomial in $q$
with rational coefficients by MacDonald's formula \cite[sec.~3.2]{M}.
This can be converted to a motive by replacing $q$ with $\ring{L}$.

By summing over the finite set of $m\in C$, we obtain the result.
\end{proof}

\section{Effective Calculations}

\begin{thm}
The virtual Chow motive $M(f,E)$ is effectively computable from the data
$f$ and $E$.
\end{thm}

\begin{rem}  \label{rem:h}
More precisely, we compute the virtual Chow motive as an
explicit linear combination of classes $[h(X)]$, where each $X$ is a
smooth projective scheme over $\ring{Q}$, $h(X)$ is the virtual Chow
motive associated with $X$, and $[h(X)]$ is its class in $\mot$. See
\cite[sec. 1.3]{DL}.  Each scheme $X$ is presented through a finite
number of coordinate patches. Each patch is given by an explicit list of
polynomials $f_1,\ldots,f_k$ in $n$ variables such that the affine patch
is isomorphic to $\ring{Q}[x_1,\ldots,x_n]/(f_1,\ldots,f_k)$. The maps
between coordinate patches are given by explicit polynomials.
\end{rem}

\begin{proof}  It is enough to prove the result when $f$ is the
characteristic function of a stable virtual set $D$ in $\lie{g}$.

Pas's algorithm is an effective procedure \cite{P}. In fact, although
Pas assumes for simplicity that the residue field has characteristic
zero, the finite set of prime characteristics that must be avoided (for
a given formula in the language) can be effectively calculated from his
algorithm.

The Presburger algorithm is an effective procedure \cite{E}.

Combining Pas's algorithm and the Presburger algorithm, we have an
effective procedure for determining effective bounds on the solutions of
each $L_i$ of Formula~\ref{eqn:vee}, for each $i\in B$.  Thus, the set
$C$ is effectively computable.  This means that an effective bound can
be obtained on the number of double cosets that must be considered.

\begin{lemma} For each $m$, the set
    $$\{(k_1,k_2,X) : \psi(m,k_1,k_2,X)\}$$
is a stable virtual subset of $G(O)\times G(O)\times \lie{g}(O)$. The
level of stability is effectively computable.
\end{lemma}

\begin{proof}  Fix $m$. For each model $L_v$, it follows directly from
compactness of $E_v$, $K_m$, and the stability of $D_v$ and $E_v$ that
the corresponding $p$-adic set is a stable subset of $G(O_v)\times
G(O_v)\times \lie{g}(O_v)$. (This is an expression of the well-known
fact that orbital integrals are locally constant on the set of regular
semisimple elements.)

Extend $\op{ord}$ to a function symbol on matrices $X=(x_{ij})$ by
defining
    $$\op{ord}(X) = \min_{ij}\op{ord}(x_{ij})$$
Let $\theta(n)=\theta_m(n)$ be the formula in Pas's language given by
    $$
    \begin{array}{lll}
    &\forall k_1\,k_2\,k_1'\,k_2'\,X\,X'.\ \\
    &\quad ( \psi(m,k_1,k_2,X)\\
        &\quad\wedge
        \op{ord}(k_1-k_1')\ge n \wedge
        \op{ord}(k_2-k_2')\ge n \wedge
        \op{ord}(X-X')\ge n)\\
    &\quad \Rightarrow \psi(m,k_1',k_2',X').
    \end{array}
    $$
The formula asserts the stability $\psi$ at level $n$. Let $\theta'(n)$
be the formula in Pas's language given by
    $$
    \theta(n)\ \wedge\ \left(\forall n'< n.\ \neg \theta(n')\right).
    $$
It asserts that $n$ is the least level for which it is stable. By the
stability of each model $L_v$ and Theorem~\ref{thm:bound}, we find that
there is a level $n$ at which almost all models are stable.  This level
is effective for the same reasons that the set $C$ above is.
\end{proof}

Now we continue with the proof that $M(f,E)$ is effectively computable.
By quantifier elimination, for each $m$, the formula $\psi(m,k_1,k_2,X)$
can be replaced by an explicit special formula (in the sense of
\cite[sec. 5.3]{DL}).   Let $n$ be the level of stability.    Truncation
$\tau_n(\psi)$ at level $n$ gives a formula in $\calL_n(G\times G\times
\lie{g})$ (the truncated arc-space of \cite[sec.~5.4]{DL}), or
equivalently the formula
    $$\psi' = \tau_n(\psi) \wedge (x\in G\times G\times \lie{g}).$$
This formula may be considered as a formula in the first-order theory of
rings.  (More explicitly, we work in the model $\ring Q[[t]]$.  We
replace quantification of $k_1$, $k_2$, $X$ with quantifiers over the
matrix coefficients, and each matrix coefficient is expanded as a
truncated power series $x_{ij} = z_0 + z_1 t + z_2 t^2 +\cdots.$)

By the definitions in \cite{DL}, we have
    $$M(f,E) = \chi_c(\psi')\ring{L}^{-(n+1)d},$$
where $d$ is the dimension of $G\times G\times\lie{g}$.

To conclude the proof, we check that $\chi_c$ can be made effective for
presented formulas $\psi'$ in the first order theory of rings.  The
procedure is as follows.  We write $\psi'$ in prenex normal form to get
a Galois formula as in \cite[sec.~2.2]{DL}.  Quantifier elimination
transforms the Galois formula into a Galois stratification.  This
procedure is effective \cite[ch.~26]{FJ}. (See also \cite[2,3,3]{DL}.)
The Galois group is presented as an explicit subgroup of a symmetric
group. Part of the data obtained from the Galois stratification is a
central function on the Galois group.

By \cite[26.12]{FJ}, Artin induction is effective, giving the central
function as an explicit rational combination of characters induced from
the trivial character on cyclic subgroups. (Fried and Jardin state the
lemma for $\ring{Q}$-valued characters, but in our case the characters
are take values in a cyclotomic extension of $\ring{Q}$.  The degree of
the cyclotomic extension is effectively computable from the degree of
the Galois cover. We can then extend the Fried-Jardin result by picking
an explicit $\ring{Q}$-basis for the cyclotomic extension.)

By refining the Galois stratification as necessary, we may assume that
the Galois stratification is affine, that is each $C/A$ is a ring cover
in the sense of \cite{FJ}.  Artin induction and the properties of
$\chi_c$ established in \cite{DL} reduce the problem to computing
$\chi_c(X/H)$, where $X$ is affine with a cyclic group action $H$, and
the data $(X,H)$ are explicitly given. The calculation of the quotient
amounts to computing the invariants of $\ring{Q}[X]$ under the finite
group $H$. There are several implementations of algorithms to compute
this ring of invariants.  (A survey of implementations appears for
example in \cite[page 2]{D}.)

To compute the class $\chi_c(X)$ where $X$ is an affine variety, embed
it in a smooth projective variety $\tilde X$.  By constructive
resolution of singularity algorithms \cite{B}, $\tilde X$ can be
effectively computed from $X$.  By the properties of the map $\chi_c$,
    $$\chi_c(X) = [h(\tilde X)] - \chi_c(D),$$
where $[h(\cdot)]$ is as in Remark~\ref{rem:h}.  The divisor $D$ has
lower dimension.  Part of the constructive resolution of singularities
gives a description of the irreducible components of $D$.  By induction
on dimension, we may assume that $\chi_c(D)$ is known.

This completes the proof that $M(f,E)$ can be effectively computed.
\end{proof}

\end{document}